\documentclass []{article}
\usepackage{amsmath,amssymb,amsthm}
\newtheorem{theorem}{Theorem}[section]
\newtheorem{lemma}[theorem]{Lemma}
\newtheorem{proposition}[theorem]{Proposition}

\theoremstyle{definition}
\newtheorem{example}[theorem]{Example}
\newtheorem{definition}[theorem]{Definition}

\begin{document}

\title{ULTRA \emph{LI}-IDEALS IN LATTICE IMPLICATION ALGEBRAS AND MTL-ALGEBRAS
\thanks{This work was
supported by the Zhejiang Provincial Natural Science Foundation of
China (Grant no. Y605389) and K.C.Wong Magna Fund in Ningbo
University.}
\author{Xiaohong Zhang, Ningbo, Keyun Qin, Chengdu, and Wieslaw A. Dudek, Wroclaw}
}
\date{}
\maketitle

\vskip 8pt

\textbf{Abstract}. A mistake concerning the ultra
\textit{LI}-ideal of a lattice implication algebra is pointed out,
and some new sufficient and necessary conditions for an
\textit{LI}-ideal to be an ultra\textit{ LI}-ideal are given.
Moreover, the notion of an \textit{LI}-ideal is extended to
MTL-algebras, the notions of a (prime, ultra, obstinate, Boolean)
\textit{LI}-ideal and an \textit{ILI}-ideal of an MTL-algebra are
introduced, some important examples are given, and the following
notions are proved to be equivalent in MTL-algebra: (1) prime
proper \textit{LI}-ideal and Boolean\textit{ LI}-ideal, (2) prime
proper \textit{LI}-ideal and \textit{ILI}-ideal, (3) proper
obstinate\textit{ LI}-ideal, (4) ultra \textit{LI}-ideal.\\

\textbf{Keywords}: lattice implication algebra, MTL-algebra,
(prime, ultra, obstinate, Boolean) \textit{LI}-ideal,\textit{
ILI}-ideal\\

\textbf{MSC2000}: 03G10, 06B10, 54E15\\

\section{Introduction}

In order to research a logical system whose propositional value is
given in a lattice, Y. Xu proposed the concept of lattice
implication algebras, and some researchers have studied their
properties and the corresponding logic systems (see \cite{Xu},
\cite{XuR}). In \cite{JunR}, Y. B. Jun et al. proposed the concept
of an \textit{LI}-ideal of a lattice implication algebra,
discussed the relationship between filters and \textit{LI}-ideals,
and studied how to generate an \textit{LI}-ideal by a set. In
\cite{Qin}, K. Y. Qin et al. introduced the notion of ultra
\textit{LI}-ideals in lattice implication algebras, and gave some
sufficient and necessary conditions for an \textit{LI}-ideal to be
ultra \textit{LI}-ideal.

The interest in the foundations of fuzzy logic has been rapidly
growing recently and several new algebras playing the role of the
structures of truth-values have been introduced. P. H\'{a}jek
introduced the system of basic logic ($BL$) axioms for the fuzzy
propositional logic and defined the class of $BL$-algebras (see
\cite{Haj}). G. J. Wang proposed a formal deductive system
$L^{\!\ast }$ for fuzzy propositional calculus, and a kind of new
algebraic structures, called $R_{0}$-algebras (see \cite{Wan1},
\cite{Wan2}). F. Esteva and L. Godo proposed a new formal
deductive system $MTL$, called the monoidal $t$-norm-based logic,
intended to cope with left-continuous $t$-norms and their
residual. The algebraic semantics for $MTL$ is based on
$MTL$-algebras (see \cite{Est}, \cite{Jen}). It is easy to verify
that a lattice implication algebra is an $MTL$-algebra. Varieties
of $MTL$-algebras are described in \cite{NEG}.

This paper is devoted to a discussion of the ultra
\textit{LI}-ideals, we correct a mistake in \cite{Qin} and give
some new equivalent conditions for an \textit{LI}-ideal to be
ultra. We also generalize the notion of an \textit{LI}-ideal to
$MTL$-algebras, introduce the notions of a (prime, ultra,
obstinate, Boolean) \textit{LI}-ideal and an \textit{ILI}-ideal of
$MTL$- algebra, give some important examples, and prove that the
following notions are equivalent in an $MTL$-algebra: (1) prime
proper \textit{LI}-ideal and Boolean \textit{ LI}-ideal, (2) prime
proper \textit{LI}-ideal and \textit{ILI}-ideal, (3) proper
obstinate \textit{ LI}-ideal, (4) ultra \textit{LI}-ideal.

\section{Preliminaries}

\begin{definition} (\cite{XuR})\label{D2.1} By a {\it lattice
implication algebra } $L$ we mean a bounded lattice
$(L,\vee,\wedge,0,1)$ with an order-reversing involution
$^{\prime}$ and a binary operation $\to$ satisfying the following
axioms:
\begin{enumerate}
\item[(I1)] \ $x\to (y\to z)=y\to (x\to z)$,
\item[(I2)] \ $x\to x=1$,
\item[(I3)] \ $x\to y=y^{\prime}\to x^{\prime}$,
\item[(I4)] \ $x\to y=y\to x=1\Longrightarrow  x=y$,
\item[(I5)] \ $(x\to y)\to y=(y\to x)\to x$,
\item[(L1)] \ $(x\vee y)\to z=(x\to z)\wedge (y\to z)$,
\item[(L2)] \ $(x\wedge y)\to z=(x\to z)\vee (y\to z)$ for all $x,y,z\in L$.
\end{enumerate}
\end{definition}

We can define a partial ordering $\le $ on a lattice implication
algebra $L$ by

\medskip \centerline{$x\le y$ if and only if $x \to y=1$.}

\medskip
For any lattice implication algebra $L$, $(L,\vee,\wedge)$ is a
distributive lattice and the De Morgan law holds, that is
\begin{enumerate}
\item[(L3)] \ $x\vee (y\wedge z)=(x\vee y)\wedge (x\vee z)$, \ \ \ \
 $x\wedge (y\vee z)=(x\wedge y)\vee (x\wedge z)$,
\item[(L4)] \ $(x\wedge y)^{\prime}=x^{\prime}\vee y^{\prime}$, \ \ \ \
 $(x\vee y)^{\prime}=x^{\prime}\wedge y^{\prime}$ for all $x,y,z\in L$.
\end{enumerate}

\begin{theorem}{\rm (\cite{XuR})}\label{T2.2} In a lattice implication
algebra $L$, the following relations hold:
\begin{enumerate}
\item[$(1)$] \ $0\to x=1$, \ \ $1\to x=x$ and $x\to 1=1$,
\item[$(2)$] \ $x$\textit{$^{\prime} $}$^{ }=x \to $0,
\item[$(3)$] \ $x\to y\le (y\to z)\to (x\to z)$,
\item[$(4)$] \ $x\vee y=(x\to y)\to y$,
\item[$(5)$] \ $x\le y$ implies $y\to z\le x\to z$ and $z\to x\le z\to
y$,
\item[$(6)$] \ $x\to (y\vee z)=(x\to y)\vee (x\to z)$,
\item[$(7)$] \ $x\to (y\wedge z)=(x\to y)\wedge (y\to z)$,
\item[$(8)$] \ $(x\to y)\vee (y\to x)=1$,
\item[$(9)$] \ $x\to (y\wedge z)=(x\to y)\wedge (x\to z)$,
\item[$(10)$] \ $x\to (y\to z)=y\to (x\to z)$,
\item[$(11)$] \ $((x\to y)\to y)\to y=x\to y$.
\end{enumerate}
\end{theorem}

From the above theorem it follows that lattice implication
algebras are strictly connected with $BCC$-algebras and
$BCK$-algebras of the form $(L,\to, 1)$ \cite{DZ1}.

For shortness, in the sequel the formula $(x\to y^{\prime}
)^{\prime}$ will be denoted by $x\otimes y $, the formula
$x^{\prime}\to y$ by $x\oplus y$.

\begin{theorem}{\rm (\cite{XuR})}\label{T2.3} In a lattice implication
algebra $L$, the relations
\begin{enumerate}
\item[$(12)$] \ $x\otimes y=y\otimes x$, \ \ $x\oplus y=y\oplus x$,
\item[$(13)$] \ $x\otimes (y\otimes z)=(x\otimes y)\otimes z$, \ \ $x\oplus (y\oplus
z)=(x\oplus y)\oplus z$,
\item[$(14)$] \ $x\otimes x^{\prime} =0$, \ \ $x\oplus
x^{\prime}=1$,
\item[$(15)$] \ $x\otimes (x\to y)=x\wedge y$,
\item[$(16)$] \ $x\to (y\to z)=(x\otimes y)\to z$,
\item[$(17)$] \ $x\le y\to z\Longleftrightarrow x\otimes y\le z$,
\item[$(18)$] \ $x\le a$ and $y\le b$ imply $x\otimes y\le a\otimes
b$ and $x\oplus y\le a\oplus b$
\end{enumerate}
hold for all $x,y,z\in L$.
\end{theorem}

\begin{definition} (\cite{JunR}) \label{D2.4} A subset $A$ of a lattice
implication algebra $L$ is called an \textit{$LI$-ideal} of $L$ if
\begin{enumerate}
\item[$(LI1)$] \ $0\in A$,
\item[$(LI2)$] \ $(x\to y)^{\prime}\in A$ and $y\in A$ imply $x\in A$
for all $x,y\in L$.
\end{enumerate}
\end{definition}

An \textit{LI}-ideal $A $ of a lattice implication algebra $L $ is
said to be {\it proper} if $A\ne L$.

\begin{theorem}\label{T2.5} {\rm (\cite{JunR}, \cite{XuR})} Let $A $ be an
\textit{LI}-ideal of a lattice implication algebra $L$, then
\begin{enumerate}
\item[{\rm (LI3)}] \ $x\in A$, $y\le x$ imply $y\in A$,
\item[{\rm (LI4)}] \ $x,y\in A$ imply $x\vee y\in A$.
\end{enumerate}
\end{theorem}

The least \textit{LI}-ideal containing a subset $A$ is called the
\textit{LI-ideal generated by} $A$ and is denoted by $\langle
A\rangle$.

\begin{theorem}\label{T2.6}{\rm (\cite{JunR}, \cite{XuR})} If $A$ is a
non-empty subset of a lattice implication algebra $L$, then

\begin{center}
$\langle A\rangle =\{x\in L\vert  a_n^{\prime}\to (\ldots\to
(a_{1}^{\prime}\to x^{\prime})\ldots )=1$ for some
$a_{1},\ldots,a_{n}\in A\}$.
\end{center}
\end{theorem}

\begin{theorem}\label{T2.7}{\rm (\cite{Qin})} Let $A$ be a subset
of a lattice implication algebra $L$. Then $A$ is an
\textit{LI}-ideal of $L$ if and only if it satisfies {\rm (LI3)}
and
\begin{enumerate}
\item[{\rm (LI5)}] \ $x\in A$ and $y\in A$ imply $x\oplus y\in A$.
\end{enumerate}
\end{theorem}

\begin{theorem}\label{T2.8}{\rm (\cite{Qin})} If $A$ is a non-empty
subset of a lattice implication algebra $L$, then

\begin{center}$\langle A\rangle =\{x\in L\vert x\le
a_{1}\oplus a_{2}\oplus\ldots\oplus a_{n}$ for some
$a_{1},...,a_{n}\in A\}$.
\end{center}
\end{theorem}

\begin{definition}\label{D2.9} (\cite{Qin}) An \textit{LI}-ideal $A$
of a lattice implication algebra $L$ is said to be {\it ultra} if
for every $x\in L$, the following equivalence holds:
\begin{enumerate}
\item[(LI6)] \ $x\in A\Longleftrightarrow x^{\prime}\notin A$.
\end{enumerate}
\end{definition}

\begin{definition}\label{D2.10} (\cite{LLX}) A non-empty subset $A$ of
a lattice implication algebra $L$ is said to be an
\textit{ILI-ideal} of $L$ if it satisfies (LI1) and
\begin{enumerate}
\item[(LI7)] \ $(((x\to y)^{\prime}\to y)^{\prime}\to z)^{\prime}\in A$ and
$z\in A$ imply $(x\to y)^{\prime}\in A$ for all $x,y,z\in L$.
\end{enumerate}
\end{definition}

\begin{theorem}\label{T2.11} {\rm (\cite{LLX})} If $A$ is an
\textit{LI}-ideal of a lattice implication algebra $L$, then the
following assertions are equivalent:
\begin{enumerate}
\item[$(i)$] \ $A$ is an \textit{ILI}-ideal of $L$,
\item[$(ii)$] \ $((x\to y)^{\prime}\to y)^{\prime}\in A$ implies $(x\to y)^{\prime}\in A$
for all $x,y,z\in L$,
\item[$(iii)$] \ $((x\to y)^{\prime}\to z)^{\prime}\in A$ implies
$((x\to z)^{\prime}\to (y\to z)^{\prime})^{\prime}\in A$ for all
$x,y,z\in L$,
\item[$(iv)$] \ $(x\to (y\to x)^{\prime})^{\prime}\in A$ implies $x\in A$ for all $x,y,z\in L$.
\end{enumerate}
\end{theorem}

\begin{definition}\label{D2.12} (\cite{Jun}) A proper
\textit{LI}-ideal $A$ of a lattice implication algebra $L$ is said
to be a {\it prime} \textit{LI}-ideal of $L$ if $x\wedge y\in A$
implies $x\in A$ or $y\in A$ for any $x,y\in L$.
\end{definition}

\begin{theorem}\label{T2.13} {\rm (\cite{LLX})} Let $A$ be a proper
\textit{LI}-ideal of a lattice implication algebra $L$. The
following assertions are equivalent:
\begin{enumerate}
\item[$(i)$] \ $A$ is a prime \textit{LI}-ideals of $L$,
\item[$(ii)$] \ $x\wedge y=0$ implies $x\in A$ or $y\in A$ for any $x,y\in L$.
\end{enumerate}
\end{theorem}

An \textit{LI}-ideal of a lattice implication algebra $L$ is
called {\it maximal}, if it is proper and not a proper subset of
any proper \textit{LI}-ideal of $L$.

\begin{theorem}\label{T2.14} {\rm (\cite{LLX})} In a lattice
implication algebra $L$, any maximal \textit{LI}-ideal must be
prime.
\end{theorem}

\begin{theorem}\label{T2.15} {\rm (\cite{LLX})} Let $L$ be a lattice
implication algebra and $A$ a proper \textit{LI}-ideal of $L$.
Then $A$ is both a prime \textit{LI}-ideal and an
\textit{ILI}-ideal of $L$ if and only if $x \in A$ or
$x$\textit{$^{\prime} $}$ \in A$ for any $x \in L$.
\end{theorem}

\begin{theorem}\label{T2.16} {\rm (\cite{LLX})} Let $L$ be a lattice
implication algebra and $A$ a proper \textit{LI}-ideal. Then $A$
is both a maximal \textit{LI}-ideal and an \textit{ILI}-ideal if
and only if for any $x,y\in L$, $x\notin A$ and $y\notin A$ imply
$(x\to y)^{\prime}\in A$ and $(y\to x)^{\prime}\in A$.
\end{theorem}

\begin{definition}\label{D2.17} (\cite{Bel}, \cite{Est}) A {\it residuated
lattice} is an algebra $(L,\wedge,\vee,\otimes,\to,0,1)$ with four
binary operations and two constants such that
\begin{enumerate}
\item[$(i)$] \ $(L,\wedge,\vee,0,1)$ is a lattice with the largest element $1$ and
the least element $0$ (with respect to the lattice ordering
$\le)$,
\item[$(ii)$] \ $(L,\otimes,1)$ is a commutative semigroup with the unit element $1$,
i.e., $\otimes$ is commutative, associative, $1\otimes x=x$ for
all $x$,
\item[$(iii)$] \ $\otimes$ and $\to$ form an adjoint pair, i.e., $z\le x\to y$
if and only if $z\otimes x\le y$ for all $x,y,z\in L$.
\end{enumerate}
\end{definition}

\begin{definition}\label{D2.18} (\cite{Est}) A residuated lattice $L$
is called an {\it $MTL$-algebra}, if it satisfies the
pre-linearity equation: $(x\to y)\vee (y\to x)=1$ for all $x,y\in
L$. An $MTL$-algebra $L$ is called an {\it $IMTL$-algebra}, if
$(a\to 0)\to 0=a$ for any $a\in L$.
\end{definition}

In the sequel $x^{\prime}$ will be reserved for $x\to 0$, $L$ for
$(L,\wedge,\vee,\otimes,\to,0,1).$

\begin{proposition}\label{P2.19} {\rm (\cite{Est}, \cite{Tur})} Let $L$ be a
residuated lattice. Then for all $x,y,z\in L$,
\begin{enumerate}
\item[{\rm (R1)}] \ $x\le y\Longleftrightarrow x\to y=1$,
\item[{\rm (R2)}] \ $x=1\to x$, \ \ $x\to (y\to x)=1$, \ \ $y\le (y\to x)\to x$,
\item[{\rm (R3)}] \ $x\le y\to z\Longleftrightarrow y\le x\to z$,
\item[{\rm (R4)}] \ $x\to (y\to z)=(x\otimes y)\to z=y\to (x\to z)$,
\item[{\rm (R5)}] \ $x\le y$ implies $z\to x\le z\to y$ and $y\to z\le
x\to z$,
\item[{\rm (R6)}] \ $z\to y\le (x\to z)\to (x\to y)$, \ \
 $z\to y\le (y\to x)\to (z \to x)$,
\item[{\rm (R7)}] \ $(x\to y)\otimes (y\to z)\le x\to z$,
\item[{\rm (R8)}] \ $x^{\prime}=x^{\prime\prime\prime}$, \ \
  $x\le x^{\prime\prime}$,
\item[{\rm (R9)}] \ $x^{\prime}\wedge y^{\prime}=(x\vee y)^{\prime}$,
\item[{\rm (R10)}] \ $x\vee x^{\prime}=1$ implies $x\wedge x^{\prime}=0$,
\item[{\rm (R11)}] \ $(\bigvee\limits_{i\in\Gamma}y_i)\to x=
 \bigwedge\limits_{i\in\Gamma }(y_i\to x)$,
\item[{\rm (R12)}] \ $x\otimes (\bigvee\limits_{i\in\Gamma }y_i)=
 \bigvee\limits_{i\in\Gamma }(x\otimes y_i)$,
\item[{\rm (R13)}] \ $x\to (\bigwedge\limits_{i\in\Gamma }y_i)=
 \bigwedge\limits_{i\in\Gamma }(x\to y_i)$,
\item[{\rm (R14)}] \ $\bigvee\limits_{i\in\Gamma }(y_i\to x)\le
 (\bigwedge\limits_{i\in\Gamma }y_i )\to x$,
 \end{enumerate}
where $\Gamma$ is a finite or infinite index set and we assume
that the corresponding infinite meets and joints exist in $L$.
\end{proposition}

\begin{proposition}\label{P2.20}{\rm (\cite{Est}, \cite{ZLi})} Let $L$ be an
$MTL$-algebra. Then for all $x,y,z\in L$,
\begin{enumerate}
\item[{\rm (M1)}] \ $x\otimes y\le x\wedge y$,
\item[{\rm (M2)}] \ $x\le y $ implies $x \otimes z\le y\otimes z$,
\item[{\rm (M3)}] \ $y\to z\le x\vee y\to x \vee z$,
\item[{\rm (M4)}] \ $x^{\prime}\vee y^{\prime}=(x\wedge
y)^{\prime}$,
\item[{\rm (M5)}] \ ($x\wedge y)\to z=(x\to z)\vee (y\to z)$,
\item[{\rm (M6)}] \ $x\vee y=((x\to y)\to y)\wedge ((y\to x)\to
x)$,
\item[{\rm (M7)}] \ $x\to (y\vee z)=(x\to y)\vee (x\to z)$,
\item[{\rm (M8)}] \ $x\wedge (y\vee z)=(x\wedge y)\vee (x\wedge z)$, \ \
 $x\vee (y\wedge z)=(x\vee y)\wedge (x\vee z)$, i.e., the lattice
structure of $L$ is distributive.
\end{enumerate}
\end{proposition}

\begin{definition}\label{D2.21} (\cite{Est}) Let $L$ be an
$MTL$-algebra. A {\it filter} is a nonempty subset $F$ of $L$ such
that
\begin{enumerate}
\item[(F1)] \ $x\otimes y\in F$ for any $x,y\in F$,
\item[(F2)] for any $x\in F$ $x\le y$ implies $y\in F$.
\end{enumerate}
\end{definition}

\begin{proposition}\label{P2.22} {\rm (\cite{Est})}
A subset $F$ of an $MTL$-algebra $L$ is a filter of $L$ if and
only if
\begin{enumerate}
\item[{\rm (F3)}] \ $1\in F$,
\item[{\rm (F4)}] \ $x\in F$ and $x\to y\in F$ imply $y\in F$.
\end{enumerate}
\end{proposition}

\section{Ultra \emph{LI}-ideals of lattice implication algebras}

In \cite{Qin}, the following result is presented: {\it let $A$ be
a subset of a lattice implication algebra $L$, then $A$ is an
ultra \textit{LI}-ideal of $L$ if and only if $A$ is a maximal
proper \textit{LI}-ideal of $L$}. The following example shows that
this result is not true.

\begin{example}\label{E3.1} Let $L=\{0,a,b,1\}$ be a set with
the Cayley table
\[
\begin{array}{|c|l|l|l|l|}\hline
 \to& 0& a& b& 1\\ \hline
   0& 1& 1& 1& 1\\ \hline
   a& b& 1& 1&1 \\ \hline
   b& a& b& 1&1 \\ \hline
   1& 0& a& b&1 \\ \hline
\end{array}
\]
For any $x\in L$, we have $x^{\prime}=x\to 0$. The operations
$\wedge$ and $ \vee $ on $L$ are defined as follows:
\[
x\vee y=(x\to y)\to y, \ \ \  x\wedge y=((x^{\prime}\to
y^{\prime})\to y^{\prime})^{\prime} .
\]
Then ($L$, $ \vee $, $ \wedge $, 0, 1) is a lattice implication
algebra. It is easy to check that $\{0\}$ is a maximal proper
\textit{LI}-ideal of $L$, but not an ultra \textit{LI}-ideal of
$L$, because $a^{\prime}=b\notin\{0\}$, but $a\notin\{0\}$.
\end{example}

Below, we give some new sufficient and necessary conditions for an
\textit{LI}-ideal to be an ultra\textit{ LI}-ideal.

\begin{theorem}\label{T3.2} Let $L$ be a lattice implication algebra and $A$ an
\textit{LI}-ideal of $L$. Then the following assertions are equivalent:
\begin{enumerate}
\item[$(i)$] \ $A$ is an ultra \textit{LI}-ideal,
\item[$(ii)$] \ $A$ is a prime proper \textit{LI}-ideal and an
\textit{ILI}-ideal,
\item[$(iii)$] \ $A$ is a proper \textit{LI}-ideal and $x\in A$ or
$x^{\prime}\in A$ for any $x\in L$,
\item[$(iv)$] \ $A$ is a maximal \textit{ILI}-ideal,
\item[$(v)$] \ $A$ is a proper \textit{LI}-ideal and for any $x,y\in L$, $x\notin A$ and
$y\notin A$ imply $(x\to y)^{\prime}\in A$ and $(y\to
x)^{\prime}\in A$.
\end{enumerate}
\end{theorem}
\begin{proof} $(i)\Longrightarrow (ii)$: $A$ is
a proper \textit{LI}-ideal, because $0\in A$, and so $1=0^{\prime}
\notin A$.

We show that $A$ is a prime \textit{LI}-ideal. Assume $x\wedge
y=0$ for some $x,y\in L$. We prove that $x\in A$ or $y\in A$. If
$x\notin A$ and $y\notin A$, then $x^{\prime}\in A$ and
$y^{\prime}\in A$, by the definition of an ultra
\textit{LI}-ideal. So, by Theorem \ref{T2.5} (LI4), we have
$x^{\prime}\vee y^{\prime}\in A$, thus $1=0^{\prime}=(x\wedge
y)^{\prime}=x^{\prime}\vee y^{\prime}\in A$. This means that
$A=L$, a contradiction. Therefore $x\wedge y$=0 implies $x\in A$
or $y\in A$. So, by Theorem \ref{T2.13}, $A$ is a prime proper
\textit{LI}-ideal.

Now we show that $A$ is an \textit{ILI}-ideal. Let $((x \to
y)^{\prime} \to y)^{\prime}\in A$. If $(x\to y)^{\prime}\notin A$,
then $x\to y\in A$ by the definition of an ultra
\textit{LI}-ideal. Since $y\le x\to y$, we have $y\in A$. From
$((x \to y)^{\prime}\to y)^{\prime}\in A$ and $y\in A$, we
conclude $(x\to y)^{\prime}\in A$, by Definition \ref{D2.4} (LI2).
This is a contradiction. Thus, $(x \to y)^{\prime}\in A$. By
Theorem \ref{T2.11} $(ii)$, $A$ is an \textit{ILI}-ideal. This
means that $(ii)$ holds.

$(ii)\Longleftrightarrow (iii)$: See Theorem \ref{T2.15}.

$(iii)\Longrightarrow (i)$: For any $x \in L$, if
$x^{\prime}\notin A$ then $x\in A$ by $(iii)$. If $x\in A$, we
prove that $x^{\prime} \notin A$. Indeed, if $x^{\prime}\in A$,
then $x \oplus x^{\prime} =1 \in A$ by Theorem \ref{T2.3} $(14)$
and Theorem \ref{T2.7} (LI5). This is a contradiction with the
fact that $A$ is a proper \textit{LI}-ideal. This means that $ A$
is an ultra \textit{LI}-ideal.

$(iv)\Longleftrightarrow (v)$: See Theorem \ref{T2.16}.

$(i)\Longrightarrow (v)$: $A$ is a proper \textit{LI}-ideal,
because $0\in A$, and so $1=0^{\prime}\notin A$.

If $x \notin A$, from $x\le y\to x$ and Theorem \ref{T2.7} (LI3),
we have $y \to x \notin A$. Thus, by the definition of an ultra
\textit{LI}-ideal, $(y\to x)^{\prime}\in A$. Similarly, from $y
\notin A $ we obtain $(x\to y)^{\prime}\in A$. That is, $(v)$
holds.

$(v)\Longrightarrow (i)$: By $(v)$, $1\notin A$. If
$x^{\prime}\notin A$, by $(v)$ we have $(1\to
x^{\prime})^{\prime}\in A$, that is $x\in A$. If $x\in A$, then
$x^{\prime}\notin A$ (see $(iii)\Longrightarrow (i)$). This means
that $A$ is an ultra \textit{LI}-ideal. The proof is complete.
\end{proof}

Remind \cite{Qin} that a subset $A$ of a lattice implication
algebra $L$ has the {\it finite additive property} if $a_{1}
\oplus a_{2}\oplus\ldots\oplus a_{n}\ne 1$ for any finite members
$a_{1},\ldots,a_{n}\in A$. $\langle A\rangle $ is a proper
\textit{LI}-ideal of $L$ if and only if $A$ has the finite
additive property.

Our theorem proves that the part results formulated in Theorem 3.7
and Corollary 3.8 in \cite{Qin} is correct. Namely we have

\begin{theorem}\label{T 3.3} If a subset $A$ of a lattice implication
algebra $L$ has the finite additive property, then there exists a
maximal \textit{LI}-ideal of $L$ containing $A$. Every proper
\textit{LI}-ideal of a lattice implication algebra can be extended
to a maximal \textit{LI}-ideal.
\end{theorem}

\section{\emph{LI}-ideals of MTL-algebras}

\begin{definition}\label{D4.1} A subset $A$ of an $MTL$-algebra $L$ is
called an \textit{LI-ideal} of $L$ if $0\in A$ and
\begin{enumerate}
\item[(LI8)] \ $(x^{\prime}\to y^{\prime})^{\prime}\in A$ and $x\in A$
imply $y\in A$ for all$x,y\in L$.
\end{enumerate}
\end{definition}

Obviously, for a lattice implication algebra $L$, (LI8) coincides
with (LI2). For a $MTL$-algebra it is not true because
$x=x^{\prime\prime}$ is not true.

An \textit{LI}-ideal $A$ of an $MTL$-algebra $L$ is said to be
{\it proper} if $A \ne L$.

\begin{lemma}\label{L4.2} {\rm (\cite{XuR} Theorem 4.1.3)}  A
non-empty subset $A$ of a lattice implication algebra $L$ is a
filter of $L$ if and only if $\,A^{\prime}=\{a^{\prime}\,|\,a\in
A\}$ is an \textit{LI}-ideal of $L$.
\end{lemma}

For $MTL$-algebras the above lemma is not true.

\begin{example}\label{E4.3} Consider the set $L=\{0,a,b,c,d,1\}$, where $0<a<b<c<d<1$, and two
operations $\otimes$, $\to$ defined by the following two tables:
\[
\begin{array}{lcr}
\begin{array}{|c|l|l|l|l|l|l|}\hline
\otimes& 0& a& b& c& d&1\\ \hline
      0& 0& 0& 0& 0&0&0\\ \hline
      a& 0& 0& 0& 0&a&a\\ \hline
      b& 0& 0& 0& b&b&b\\ \hline
      c& 0& 0& b& c&c&c\\ \hline
      d& 0& a& b& c&d&d\\ \hline
      1& 0& a& b& c&d&1\\ \hline
\end{array}
&\rule{20mm}{0mm}&
\begin{array}{|c|l|l|l|l|l|l|}\hline
\to& 0& a& b& c&d &1\\ \hline
  0& 1& 1& 1& 1&1 &1\\ \hline
  a& c& 1& 1& 1&1 &1\\ \hline
  b& b& b& 1& 1&1 &1\\ \hline
  c& a& a& b& 1&1 &1\\ \hline
  d& 0& a& b& c&1 &1\\ \hline
  1& 0& a& b& c&d &1\\ \hline
\end{array}
\end{array}
\]
If we define on $L$ the operations $\wedge$ and $\vee$ as $\min$
and $\max$, respectively, then $(L,\wedge,\vee,\otimes,\to,0,1)$
will be an $MTL$-algebra. Obviously, $A=\{0,a,b,c,d,1\}$ is a
filter of $L$, but $A^{\prime}=\{0,a,b,c,1\}$ is not an
\textit{LI}-ideal of $L$, since
\[
(0^{\prime}\to d^{\prime})^{\prime}=1\in A \ \ {\rm and } \ \ 0\in
A, \ \ {\rm but } \ \ d\notin A.
\]
Moreover, $B=\{1,c\}$ is not a filter of $L$, because $c\to d=1\in
B$ and $c\in B$, $d\notin B$. By the following MATHEMATICA
program, we can verify that $B^{\prime}=\{0,a\}$ is an \textit{LI}-ideal of $L$:\\

M1={\{}{\{}6,6,6,6,6,6{\}},{\{}4,6,6,6,6,6{\}},{\{}3,3,6,6,6,6{\}},{\{}2,2,3,6,6,6{\}},{\{}1,2,3,4,6,6{\}},{\{}1,2,3,4,5,6{\}}{\}};

a1=0;

For[i=1, i$<$7, i++, For[j=1, j$<$7, j++,

$\ \ \ \ $If[(i==1$\vert \vert $i==2) {\&}{\&}
(M1[[M1[[M1[[i,1]],M1[[j,1]]]],1]]==1$\vert \vert $

$\ \ \ \ \ \ \ \ $M1[[M1[[M1[[i,1]],M1[[j,1]]]],1]]==2) {\&}{\&}
(j!=1{\&}{\&} j!=2), a1++]]];

If[a1==0, Print["true"], Print["false"]]
\end{example}

From Example \ref{E4.3} we see that \textit{LI}-ideals have a
proper meaning in $MTL$-algebras.

\begin{theorem}\label{T4.4} Let $A$ be an \textit{LI}-ideal of an $MTL$-algebra $L$, then
\begin{enumerate}
\item[{\rm (LI3)}] \ if $x\in A$, $y\le x$, then $y \in A$,
\item[{\rm (LI9)}] \ if $x\in A$, then $x^{\prime \prime}\in A$,
\item[{\rm (LI4)}] \ if $x,y\in A$, then $x \vee y \in A$.
\end{enumerate}
\end{theorem}
\begin{proof} Assume $x\in A$, $y\le x$. From $y\le x$, by
Proposition \ref{P2.19} (R5), we have $x\to 0\le y\to 0$, i.e.,
$x^{\prime}\le y^{\prime}$. By Proposition \ref{P2.19} (R1),
$x^{\prime}\to y^{\prime}=1$. Then $(x^{\prime}\to
y^{\prime})^{\prime}=1^{\prime}=0\in A$ and $x\in A$, and by (LI8)
we get $y \in A$. This means that (LI3) holds.

Suppose $x\in A$. By Proposition \ref{P2.19} (R8) we have
$(x^{\prime}\to (x^{\prime\prime})^{\prime})^{\prime}=
(x^{\prime}\to x^{\prime})^{\prime}=1^{\prime}=0\in A$. Applying
(LI8) we get $x^{\prime\prime}\in A$, i.e., (LI9) holds.

Assume $x,y\in A$. By Proposition \ref{P2.19} (R2) we have
$y^{\prime}\le x^{\prime}\to y^{\prime}$. So, $(x^{\prime}\to
y^{\prime})^{\prime}\le y^{\prime\prime}$ by (R5). Whence, by
$y\in A$ and (LI9), we obtain $y^{\prime\prime}\in A$. From this
and $(x^{\prime}\to y^{\prime})^{\prime}\le y^{\prime\prime}$,
using (LI3) we get $(x^{\prime}\to y^{\prime})^{\prime}\in A$.
Thus

\medskip
\begin{tabular}{rlcr}
$(x^{\prime}\to (x\vee
y)^{\prime})^{\prime}\!\!\!\!\!\!$&$=(x^{\prime} \to
(x^{\prime}\wedge y^{\prime}))^{\prime}$&$\rule{40mm}{0mm}$&(by
(R9))\\[4pt]
&$=((x^{\prime}\to x^{\prime})\wedge (x^{\prime}\to
y^{\prime}))^{\prime}$&& (by (R13))\\[4pt]
&$=(1\wedge (x^{\prime}\to y^{\prime}))^{\prime}$&& (by (R1))\\[4pt]
&$=(x^{\prime}\to y^{\prime})^{\prime}\in A$.
\end{tabular}

\medskip
\noindent From this and $x\in A$, using (LI8), we deduce $x\vee
y\in A$, i.e., (LI4) holds.

The proof is complete.
 \end{proof}

\begin{definition}\label{D4.5} An \textit{LI}-ideal $A$ of an
$MTL$-algebra $L$ is said to be an \textit{ILI-ideal} of $L$ if it
satisfies
\begin{enumerate}
\item[(LI10)] \ $(x\to (y\to x)^{\prime})^{\prime}\in A$ implies $x\in A
$ for all $x,y,z\in L$.
\end{enumerate}
\end{definition}

\begin{example}\label{E4.6} Let $L=\{0,a,b,1\}$, where $0<a<b<1$, be a set with
the Cayley tables:
\[
\begin{array}{lcr}
\begin{array}{|c|l|l|l|l|}\hline
\otimes& 0& a& b& 1\\ \hline
      0& 0& 0& 0& 0\\ \hline
      a& 0& a& a& a\\ \hline
      b& 0& a& a& b\\ \hline
      1& 0& a& b& 1\\ \hline
\end{array}
&\rule{20mm}{0mm}&
\begin{array}{|c|l|l|l|l|}\hline
\to & 0& a& b& 1\\ \hline
   0& 1& 1& 1& 1\\ \hline
   a& 0& 1& 1& 1\\ \hline
   b& 0& b& 1& 1\\ \hline
   1& 0& a& b& 1\\ \hline
\end{array}
\end{array}
\]
Defining the operations $\wedge $ and $\vee$ on $L$ as $\min$ and
$\max$, respectively, we obtain an $MTL$-algebra $(L,\wedge,\vee,
\otimes,\to ,0,1)$ in which $A=\{0\}$ is an \textit{ILI}-ideal of
$L$.
\end{example}

In Example \ref{E4.3}, $\{0,a\}$ is an \textit{LI}-ideal, but it
is not an \textit{ILI}-ideal of $L$, because
\[
(b\to (1\to b)^{\prime})^{\prime}=0\in\{0,a\}, \ \ {\rm  but } \ \
b\notin\{0,a\}.
\]

\begin{theorem}\label{T4.7} For each an $I$\textit{LI}-ideal $A$ of an $MTL$-algebra
$L$ we have
\begin{enumerate}
\item[{\rm (LI11)}] \ $x\wedge x^{\prime}\in A$ for each $x\in L$.
\end{enumerate}
\end{theorem}
\begin{proof} Indeed, for all $x\in L$ we get

\medskip\hspace*{-5mm}
\begin{tabular}{rllll}
$((x\wedge x^{\prime})\to (1\to (x\wedge
x^{\prime}))^{\prime})^{\prime}\!\!\!\!\!\!$&$=((x\wedge
x^{\prime})\to
(x\wedge x^{\prime})^{\prime})^{\prime}$&&\\[4pt]
&$=((x\wedge x^{\prime})\to (x^{\prime}\vee
x^{\prime\prime}))^{\prime}$&(by Proposition \ref{P2.20}
(M4))\\[4pt]
&$=(((x\wedge x^{\prime})\to x^{\prime})\vee ((x\wedge x^{\prime})
\to x^{\prime\prime}))^{\prime}$&(by Proposition \ref{P2.20}
(M7))\\[4pt]
&$=(1\vee ((x\wedge x^{\prime})\to x^{\prime\prime}))^{\prime}$&
(by Proposition \ref{P2.19} (R1))\\[4pt]
&$=1^{\prime}=0\in A$.
\end{tabular}

\medskip\noindent
From this, applying (LI10), we deduce (LI11).
 \end{proof}

\begin{definition}\label{D4.8} An \textit{LI}-ideal $A$ satisfying (LI11) is called
{\it Boolean}.
\end{definition}

\begin{theorem}\label{T4.9} If $A$ is a Boolean \textit{LI}-ideal of an $MTL$-algebra $L$, then
\begin{enumerate}
\item[{\rm (LI12)}] \ $(x\to x^{\prime})^{\prime}\in A$ implies $x\in A$.
\end{enumerate}
\end{theorem}
\begin{proof} According to the assumption $x\wedge x^{\prime}\in A$ for all $x\in L$.
Let $(x\to x^{\prime})^{\prime}\in A$. Then

\medskip
\begin{tabular}{rlrr}
$((x\wedge x^{\prime})^{\prime}\to
x^{\prime})^{\prime}\!\!\!\!\!\!$
 &$=(x\to (x\wedge x^{\prime})^{\prime\prime})^{\prime}$&&(by Proposition \ref{P2.19}
(R4))\\[4pt]
&$=(x\to (x^{\prime\prime}\wedge
x^{\prime\prime\prime}))^{\prime}$&&(by Propositions \ref{P2.19}
(R9) and \ref{P2.20} (M4))\\[4pt]
&$=((x\to x^{\prime\prime})\wedge (x\to
x^{\prime\prime\prime}))^{\prime}$&&(by Proposition \ref{P2.19}
(R13))\\[4pt]
&$=(1\wedge (x\to x^{\prime}))^{\prime}$&&(by Proposition
\ref{P2.19}
(R8), (R1))\\[4pt]
&$=(x\to x^{\prime})^{\prime}\in A$.
\end{tabular}

\medskip\noindent
Now, applying (LI8) we get $x\in A$, which completes the proof.
\end{proof}

\begin{theorem}\label{T4.10} For \textit{LI}-ideals of $MTL$-algebras the conditions
{\rm (LI10)} are equivalent {\rm (LI11)}.
\end{theorem}
\begin{proof} (LI10)$\Longrightarrow$(LI11): See Theorem \ref{T4.7}.

(LI11)$\Longrightarrow$ (LI10): Let $(x\to (y\to
x)^{\prime})^{\prime}\in A$. Then

\medskip\hspace*{-8mm}
\begin{tabular}{rlll}
$((x\to (y\to x)^{\prime})^{\prime\prime}\to (x\to x^{\prime}
)^{\prime\prime})^{\prime}\!\!\!\!\!\!$&$=((x\to
x^{\prime})^{\prime}\to (x \to (y \to
x)^{\prime})^{\prime})^{\prime}$&(by Proposition
\ref{P2.19} (R4), (R8))\\[4pt]
&$\le((x\to (y\to x)^{\prime})\to (x\to x^{\prime}))^{\prime}$&
(by Proposition \ref{P2.19} (R6))\\[4pt]
&$\le((y\to x)^{\prime}\to x^{\prime})^{\prime}$&(by Proposition
\ref{P2.19} (R6))\\[4pt]
&$\le (x\to (y\to x))^{\prime}$&(by Proposition \ref{P2.19} (R6))\\[4pt]
&$=1^{\prime}=0 \in A$.&(by Proposition \ref{P2.19} (R2))
\end{tabular}

\medskip\noindent
This, by (LI8), implies ($x\to x^{\prime})^{\prime}\in A$, whence,
using (LI12), we obtain $x\in A$. So, (LI11) implies (LI10).
\end{proof}

\begin{definition}\label{D4.11} A proper \textit{LI}-ideal $A$ of an
$MTL$-algebra $L$ is said to be a {it prime} if $x\wedge y\in A$
implies $x\in A$ or $y\in A$ for any $x,y\in L$.
\end{definition}

\begin{theorem}\label{T 4.12} A proper \textit{LI}-ideal $A$ of a
$MTL$-algebra $L$ is prime if and only if for all $x,y\in L$ we
have $(x\to y)^{\prime}\in A$ or $(y\to x)^{\prime}\in A$.
\end{theorem}
\begin{proof} Assume that an \textit{LI}-ideal $A$ of $L$ is prime. Since
\[(x\to y)^{\prime}\wedge (y\to
x)^{\prime}=((x\to y)\vee (y\to x))^{\prime}=1^{\prime}=0\in A
\]
for all $x,y\in L$, the assumption on $A$ implies $(x\to
y)^{\prime}\in A$ or $(y \to x)^{\prime}\in A$.

Conversely, let $A$ be a proper \textit{LI}-ideal of $L$ and let
$x\wedge y\in A$. Assume that $(x\to y)^{\prime}\in A$ or $(y\to
x)^{\prime}\in A$ for $x,y\in L$. If $(x\to y)^{\prime}\in A$,
then

\medskip
\begin{tabular}{rlrr}
$((x\wedge y)^{\prime}\to
x^{\prime})^{\prime}\!\!\!\!\!\!$&$=((x^{\prime}\vee
y^{\prime})\to
x^{\prime})^{\prime}$&$\rule{20mm}{0mm}$&(by Proposition \ref{P2.20} (M4))\\[4pt]
&$=((x^{\prime}\to x^{\prime})\wedge (y^{\prime}\to x^{\prime}
))^{\prime}$&&(by Proposition \ref{P2.19} (R11))\\[4pt]
&$=(1\wedge (y^{\prime}\to x^{\prime}))^{\prime}$&&(by Proposition
\ref{P2.19} (R1))\\[4pt]
&$=(y^{\prime}\to x^{\prime})^{\prime}\le (x\to y)^{\prime}\in
A$.&&(by Proposition \ref{P2.19} (R6))
\end{tabular}

\medskip\noindent
So, $((x \wedge y)^{\prime}\to x^{\prime})^{\prime}\in A$ (Theorem
\ref{T4.4} (LI3)), which together with $x\wedge y\in A$ and the
definition of an $LI$-ideal, gives $x \in A$.

Similarly, from $(y\to x)^{\prime}\in A$ we can obtain $y\in A$.

This means that $A$ is a prime \textit{LI}-ideal of $L$. The proof
is complete.
 \end{proof}

\begin{theorem}\label{T4.13} Let $A$ be an \textit{LI}-ideal of an
$MTL$-algebra $L$. Then $A$ is both a prime \textit{LI}-ideal and
a Boolean \textit{LI}-ideal of $L$ if and only if $x\in A$ or
$x^{\prime}\in A$ for any $x\in L$.
\end{theorem}
\begin{proof} Assume that for all $x\in L$ we have $x\in A$ or $x^{\prime}\in A$.
At first we show that an $LI$-ideal $A$ is prime. For this let
$x\wedge y\in A$. If $x\notin A$, then $x^{\prime}\in A$. Hence

\medskip
\begin{tabular}{rlrrr}
$((x\wedge y)^{\prime}\to y^{\prime})^{\prime}\!\!\!\!\!\!$
&$=((x^{\prime}\vee y^{\prime})\to
y^{\prime})^{\prime}$&$\rule{25mm}{0mm}$&(by
Proposition \ref{P2.20} (M4))\\[4pt]
&$=((x^{\prime}\to y^{\prime})\wedge (y^{\prime}\to
y^{\prime}))^{\prime}$&&(by Proposition \ref{P2.19} (R11))\\[4pt]
&$=((x^{\prime}\to y^{\prime})\wedge 1)^{\prime}$&&(by Proposition
\ref{P2.19} (R1))\\[4pt]
&$=(x^{\prime}\to y^{\prime})^{\prime}\le (y\to x)^{\prime}$&& (by
Proposition \ref{P2.19} (R6))\\[4pt]
&$\le x^{\prime}\in A$.&&(by Proposition \ref{P2.19} (R2))
\end{tabular}

\medskip\noindent
So, $((x\wedge y)^{\prime}\to y^{\prime})^{\prime}\in A$, by
Theorem \ref{T4.4} (LI3). From this, $x\wedge y\in A$ and
Definition \ref{D4.1} we get $y\in A$. This proves that an
$LI$-ideal $A$ is prime. To prove that it is Boolean observe that
$x\wedge x^{\prime}\le x^{\prime}$ implies $x\wedge x^{\prime}\le
x$, whence, by Theorem \ref{T4.4} (LI3), we obtain $x\wedge
x^{\prime}\in A$. Thus $A$ is Boolean.

Conversely, is an $LI$-ideal $A$ is both prime and Boolean, then
by Definition \ref{D4.8}, for all $x\in L$ we have $x\wedge
x^{\prime}\in A$. Hence $x\in A$ or $x^{\prime}\in A$, by
Definition \ref{D4.11}. This completes the proof.
 \end{proof}

\begin{definition}\label{D4.14} An \textit{LI}-ideal $A$ of an
$MTL$-algebra $L$ is said to be {\it ultra} if for every $x \in L$
\begin{enumerate}
\item[(LI6)] \ $x\in A\Longleftrightarrow x^{\prime}\notin A$.
\end{enumerate}
\end{definition}

It is easy to verify the following proposition is true.

\begin{proposition}\label{P4.15} Each ultra \textit{LI}-ideal of an $MTL$-algebra is a proper
\textit{LI}-ideal.
\end{proposition}

\begin{definition}\label{D4.16} An \textit{LI}-ideal $A$ of an
$MTL$-algebra $L$ is said to be {\it obstinate} if for all $x,y\in
L$
\begin{enumerate}
\item[(LI13)] \ $x\notin A$ and $y \notin A$
imply $(x\to y)^{\prime}\in A$ and $(y\to x)^{\prime}\in A$.
\end{enumerate}
\end{definition}

\begin{theorem}\label{T4.17} For an \textit{LI}-ideal $A$ of an $MTL$-algebra $L$ the
following conditions are equivalent:
\begin{enumerate}
\item[$(i)$] \ $A$ is an ultra \textit{LI}-ideal,
\item[$(ii)$] \ $ A$ is a proper \textit{LI}-ideal and for any $x\in L$, $x\in A$ or
$x^{\prime}\in A$,
\item[$(iii)$] \ $A$ is a prime proper \textit{LI}-ideal and a Boolean
\textit{LI}-ideal,
\item[$(iv)$] \ $A$ is a prime proper \textit{LI}-ideal and an
\textit{ILI}-ideal,
\item[$(v)$] \ $ A$ is a proper obstinate\textit{ LI}-ideal.
\end{enumerate}
\end{theorem}
\begin{proof} $(i)\longrightarrow (ii)$: Obvious.

$(ii)\longrightarrow (i)$: If $x^{\prime}\notin A$, then $x\in A$,
by $(ii)$. Similarly, if $x\in A$, that must be $x^{\prime}\notin
A$. If not, i.e., if $x^{\prime}\in A$, then, by Proposition
\ref{P2.19} (R8), we have
\[
(x^{\prime}\to 1^{\prime})^{\prime}=(x^{\prime}\to 0)^{\prime}
=x^{\prime\prime\prime}=x^{\prime}\in A,
\]
which together with $x\in A$ and (LI8) implies $1\in A$. This, by
Theorem \ref{T4.4} (LI3), gives $A=L$. This is a contradiction,
because an \textit{LI}-ideal $A$ is proper. Obtained contradiction
proves that $x\in A$ implies $x^{\prime}\notin A$. So, $A$ is an
ultra \textit{LI}-ideal.

$(ii)\Longleftrightarrow (iii)$: See Theorem \ref{T4.13}.

$(iv)\Longrightarrow (iii)$: See Theorem \ref{T4.7}.

$(iii)\Longrightarrow (iv)$: See Theorem \ref{T4.10}.

$(v)\Longrightarrow (ii)$: Since $A$ is a proper
\textit{LI}-ideal, $1\notin A$. If $x\notin A$, then $(1\to
x)^{\prime}=x^{\prime}\in A$, by Definition \ref{D4.16}. This
means that $(ii)$ holds.
\newpage
$(ii)\Longrightarrow (v)$: It suffices to show that $A$ is
obstinate. Let $x\notin A$ and $y \notin A$. Then, according to
$(ii)$, we have $x^{\prime}\in A$ and $y^{\prime}\in A$. Thus

\medskip
\begin{tabular}{rlrr}
$(y^{\prime\prime}\to (x\to
y)^{\prime\prime})^{\prime}\!\!\!\!\!\!$
 &$=((x\to y)^{\prime}\to
y^{\prime\prime\prime})^{\prime}$&$\rule{15mm}{0mm}$ &
 (by Proposition \ref{P2.19} (R4))\\[4pt]
&$=((x\to y)^{\prime}\to y^{\prime})^{\prime}$&&(by Proposition
\ref{P2.19} (R8))\\[4pt]
&$=(y\to (x\to y)^{\prime\prime})^{\prime}$&&(by Proposition
\ref{P2.19} (R4))\\[4pt]
&$\le (y\to (x\to y))^{\prime}$&&(by (R8), $x\to y\le (x\to
y)^{\prime\prime}$ and (R5))\\[4pt]
&$=1^{\prime}=0\in A$.&&(by Proposition \ref{P2.19} (R2))
\end{tabular}

\medskip\noindent
This together with $y^{\prime}\in A$ and Definition \ref{D4.1}
implies $(x\to y)^{\prime}\in A$.

Similarly, we obtain $(y\to x)^{\prime}\in A$. So, $A$ is
obstinate.

The proof is complete.
\end{proof}

Authors' addresses: X. H. Zhang, Department of Mathematics,
Faculty of Science, Ningbo University, Ningbo 315211, Zhejiang
Province, P. R. China, e-mail: zxhonghz@263.net; K. Y. Qin,
Department of Applied Mathematics, Southwest Jiaotong University,
Chengdu, Sichuan 610031, P. R. China, e-mail: keyunqin@263.net;
W.A.Dudek, Institute of Mathematics and Computer Science, Wroclaw
University of Technology, Wybrzeze Wyspianskiego 27, 50-370
Wroclaw, POLAND, e-mail: dudek@im.pwr.wroc.pl


\begin{thebibliography}{}
\bibitem{Bel} \emph{R.Belohlavek}: Some properties of residuated lattices.
Czechoslovak Math. J. {\it 53(128)} (2003), $161-171.$ Zbl
1014.03510
\bibitem{DZ1} \emph{W.A.Dudek, X.H.Zhang}: On ideals and congruences in
$BCC$-algebras. Czechoslovak Math. J. {\it 48(123)} (1998),
$21-29.$ Zbl 0927.06013
\bibitem{Est} \emph{F.Esteva, L.Godo}: Monoidal $t$-norm based
logic: Towards a logic for left-continuous $t$-norms. Fuzzy Sets
and Systems {\it 124} (2001), $271-288.$ Zbl 0994.03017
\bibitem{Haj} \emph{P.H\'ajek}: Metamathematics of Fuzzy Logic. Kluwer
Academic Publishers, 1998. Zbl 0937.03030
\bibitem{Jen} \emph{S.Jenei, F.Montagna}: A proof of standard completeness for
Esteva and Godo's logic $MTL$. Studia Logica {\it 70} (2002),
$183-192.$ Zbl 0997.03027
\bibitem{Jun} \emph{Y.B.Jun}: On $LI$-ideals and prime $LI$-ideals of lattice
implication algebras. J. Korean Math. Soc. {\it 36} (1999),
$369-380.$ Zbl 0919.03050
\bibitem{JunR} \emph{Y.B.Jun, E.H.Roh, Y.Xu}: $LI$-ideals in lattice
implication algebras. Bull. Korean Math. Soc. {\it 35} (1998),
$13-24.$ Zbl 0903.03037
\bibitem{JunX} \emph{Y.B.Jun, Y.Xu}: Fuzzy $LI$-ideals in lattice
implication algebras. J. Fuzzy Math. {\it 7} (1999), $997-1003.$
Zbl 0972.03550
\bibitem{LLX} \emph{Y.L.Liu, S.Y.Liu, Y.Xu, K.Y.Qin}: $ILI$-ideals and prime
$LI$-ideals in lattice implication algebras. Information Sciences
{\it 155} (2003), $157-175.$ Zbl 1040.03046
\bibitem{NEG} \emph{C.Noguera, F.Esteva, J.Gispert}: On some varieties of $MTL$-algebras. Log. J. IGPL
{\it 13} (2005), $443-466.$ Zbl 1078.03051
\bibitem{Qin} \emph{K.Y.Qin, Y.Xu, Y.B.Jun}: Ultra $LI$-ideals in lattice
implication algebras. Czechoslovak Math. J. {\it 52 (127)} (2002),
$463-468.$ Zbl 1012.03061
\bibitem{Tur} \emph{E.Turunen}: Boolean deductive systems of $BL$-algebras. Arch.
Math. Logic {\it 40} (2001), $467-473.$ Zbl 1030.03048
\bibitem{Wan1} G.J.Wang: Non-classical Mathematical Logic and Approximate
Reasoning, (Chinese), Beijing, Science Press, 2000.
\bibitem{Wan2} \emph{G.J.Wang}: $MV$-algebras, $BL$-algebras, $R_0$-algebras and
multiple-valued logic. (Chinese), Fuzzy Systems and Mathematics
{\it 16} (2002), No.2, $1-15.$ MR 1911031
\bibitem{Xu} \emph{Y.Xu}: Lattice implication algebras. J. South West Jiaotong
University. {\it 1} (1993), $20-27.$
\bibitem{XuQ} \emph{Y.Xu, K.Y.Qin}: On filters of lattice implication
algebras. J. Fuzzy Math. {\it 1} (1993), $251-260.$ Zbl 0787.06009
\bibitem{XuR} \emph{Y.Xu, D.Ruan, K.Y.Qin, J.Li}: Lattice-valued Logic. An alternative
approach to treat fuzziness and incomparability. Studies in
Fuzzines and soft Computing 132, Springer. 2003. Zbl 1048.03003
\bibitem{ZLi} \emph{X.H.Zhang, W.H.Li}: On fuzzy logic algebraic system $MTL$.
Advances in Systems and Applications {\it 5} (2005), $475-483.$
\end{thebibliography}
\end{document}